\documentclass[leqno,12pt]{amsart}
\usepackage{amssymb} 
\theoremstyle{plain} 
\newtheorem{theorem}{\indent\sc Theorem}[section] 

\newtheorem{proposition}[theorem]{\indent\sc Proposition}

\theoremstyle{definition} 

\begin{document}

\title[Coupled Painlev\'e systems]{Coupled Painlev\'e systems in dimension four with affine Weyl group symmetry of types $A_4^{(2)}$ and $A_1^{(1)}$ \\}

\renewcommand{\thefootnote}{\fnsymbol{footnote}}
\footnote[0]{2000\textit{ Mathematics Subjet Classification}.
34M55; 34M45; 58F05; 32S65.}

\maketitle

\begin{abstract}
We find a two-parameter family of coupled Painlev\'e systems in dimension four with affine Weyl group symmetry of type $A_4^{(2)}$.  For a degenerate system of $A_4^{(2)}$ system, we also find a one-parameter family of coupled Painlev\'e systems in dimension four with affine Weyl group symmetry of type $A_1^{(1)}$. We show that for each system, we give its symmetry and holomorphy conditions. These symmetries, holomorphy conditions and invariant divisors are new. Moreover, we find a one-parameter family of partial differential systems in three variables with $W(A_1^{(1)})$-symmetry. We show the relation between its polynomial Hamiltonian system and an autonomous version of the system of type $A_1^{(1)}$.
\end{abstract}

\section{Introduction}

In \cite{Sasa11}, we find a four-parameter (resp. three-parameter) family of ordinary differential systems in dimension four with affine Weyl group symmetry of type $A_7^{(2)}$ (resp. $A_5^{(2)}$). These systems are equivalent to the polynomial Hamiltonian systems, and can be considered to be 2-coupled Painlev\'e systems in dimension four.

We will complete the study of the below problem in a series of papers, for which this paper is the third, resulting in a series of equations for the remaining affine root systems of type $A_{2l+2}^{(2)}$.

{\bf Problem} \quad For each affine root system $X_i^{(2)}$ with affine Weyl group $W(X_i^{(2)})$, find a system of differential equations for which $W(X_i^{(2)})$ acts as its B{\"a}cklund transformations.

This paper is the stage in this project where we find a 2-parameter family of coupled Painlev\'e systems in dimension four with affine Weyl group symmetry of type $A_4^{(2)}$ given by
\begin{align}\label{1}
\begin{split}
\frac{dx}{dt}&=\frac{\partial H}{\partial y}, \quad \frac{dy}{dt}=-\frac{\partial H}{\partial x}, \quad \frac{dz}{dt}=\frac{\partial H}{\partial w}, \quad \frac{dw}{dt}=-\frac{\partial H}{\partial z}
\end{split}
\end{align}
with the polynomial Hamiltonian
\begin{align}\label{2}
\begin{split}
H =&2H_{II}(x,y,t,\alpha_1)+{H_{II}}^{auto}(z,w,t,\alpha_0)+xw+2yzw\\
=&2xy^2+2x^2+2tx-2\alpha_1 y+z^2 w-\frac{w^2}{2}+\alpha_0 z+xw+2yzw.
\end{split}
\end{align}
Here $x,y,z$ and $w$ denote unknown complex variables, and $\alpha_0,\alpha_1,\alpha_2$ are complex parameters satisfying the relation:
\begin{equation}\label{3}
\alpha_0+2\alpha_1+2\alpha_2=1.
\end{equation}
The symbol $H_{II}$ denotes the second Painlev\'e Hamiltonian given by
\begin{align}
&H_{II}(x,y,t,\alpha_1)=xy^2+x^2+tx-\alpha_1 y,
\end{align}
and the symbol ${H_{II}}^{auto}$ denotes the autonomous version of the second Painlev\'e system given by
\begin{align}
&{H_{II}}^{auto}(z,w,t,\alpha_0)=z^2 w-\frac{w^2}{2}+\alpha_0 z.
\end{align}
Of course, the Hamiltonian itself is the first integral.

We remark that for this system we tried to seek its first integrals of polynomial type with respect to $x,y,z,w$. However, we can not find. Of course, the Hamiltonian $H$ is not its first integral.

We also remark that the system \eqref{1} can be obtained by connecting the pair of the invariant divisors $(x+y+w+t,y)$ and $(x-z,z)$ for the canonical variables $(x,y,z,w)$ (see figure 1) in the system of type $A_4^{(1)}$ (see section 5 in \cite{Sasa1}).

This is the second example which gave higher-order Painlev\'e equations of type $A_{4}^{(2)}$. 

{\bf Problem} \quad It is still an open question whether the system \eqref{1} is equivalent to Ramani's equation of type $A_4^{(2)}$.

 For a degenerate system of $A_4^{(2)}$ system, we also find a one-parameter family of coupled Painlev\'e systems in dimension four with affine Weyl group symmetry of type $A_1^{(1)}$. 

We show that for each system, we give its symmetry and holomorphy conditions. These symmetries, holomorphy conditions and invariant divisors are new.

Finally, we find a 1-parameter family of partial differential systems in three variables given by (see Section 4)
\begin{equation}
  \left\{
  \begin{aligned}
   dq_1 =&\frac{\partial K_1}{\partial p_1}dt_1+\frac{\partial K_2}{\partial p_1}dt_2+\frac{\partial K_3}{\partial p_1}dt_3,\\
   dp_1 =&-\frac{\partial K_1}{\partial q_1}dt_1-\frac{\partial K_2}{\partial q_1}dt_2-\frac{\partial K_3}{\partial q_1}dt_3,\\
   dq_2 =&\frac{\partial K_1}{\partial p_2}dt_1+\frac{\partial K_2}{\partial p_2}dt_2+\frac{\partial K_3}{\partial p_2}dt_3,\\
   dp_2 =&-\frac{\partial K_1}{\partial q_2}dt_1-\frac{\partial K_2}{\partial q_2}dt_2-\frac{\partial K_3}{\partial q_2}dt_3
   \end{aligned}
  \right. 
\end{equation}
with the polynomial Hamiltonians
\begin{align}
\begin{split}
K_1=&q_1p_1^2+q_1^2-\alpha_0 p_1+\frac{1}{4}q_2^2-\frac{1}{4}p_2^2+p_1q_2p_2,\\
K_2=&q_2^2p_2^2-\frac{1}{4}q_2^2+\frac{1}{4}p_2^2-2\alpha_0 q_2p_2+q_1q_2^2+q_1p_2^2-p_1q_2p_2+p_1^2q_2^2,\\
K_3=&\frac{1}{2}q_1^2p_1^4+q_1^3p_1^2+\frac{1}{2}q_1^4-\alpha_0 q_1p_1^3+\alpha_1 q_1^2p_1+\frac{1}{2}\alpha_0^2 p_1^2+\frac{1}{32}q_2^4+\frac{1}{32}p_2^4-\frac{1}{16}q_2^2p_2^2+q_1p_1^3q_2p_2\\
&+\frac{1}{2}p_1^2q_2^2p_2^2-\frac{1}{4}q_1p_1^2p_2^2-\frac{1}{4}p_1q_2p_2^3+q_1^2p_1q_2p_2+\frac{1}{4}q_1p_1^2q_2^2+\frac{1}{4}p_1q_2^3p_2-\alpha_0 p_1^2q_2p_2\\
&+\frac{1}{4}\alpha_0 p_1p_2^2+\frac{1}{4}\alpha_1 p_1q_2^2.
\end{split}
\end{align}
Here, the parameters $\alpha_i$ satisfy the relation $\alpha_0+\alpha_1=0$. This system is new.

We see that this system satisfies the compatibility conditions$:$
\begin{equation}
\frac{\partial }{\partial t_j} \frac{\partial q_1}{\partial t_i}=\frac{\partial }{\partial t_i} \frac{\partial q_1}{\partial t_j}, \quad \frac{\partial }{\partial t_j} \frac{\partial p_1}{\partial t_i}=\frac{\partial }{\partial t_i} \frac{\partial p_1}{\partial t_j}, \quad \frac{\partial }{\partial t_j} \frac{\partial q_2}{\partial t_i}=\frac{\partial }{\partial t_i} \frac{\partial q_2}{\partial t_j}, \quad \frac{\partial }{\partial t_j} \frac{\partial p_2}{\partial t_i}=\frac{\partial }{\partial t_i} \frac{\partial p_2}{\partial t_j},
\end{equation}
where $i,j=1,2,3$, and this system has $K_1,K_2$ and $K_3$ as its first integrals, and three Hamiltonians $K_1,K_2$ and $K_3$ satisfy
\begin{equation}
\{K_1,K_2\}=\{K_1,K_3\}=\{K_2,K_3\}=0,
\end{equation}
where
\begin{equation}
\{K_i,K_j\}=\frac{\partial K_i}{\partial p_1}\frac{\partial K_j}{\partial q_1}-\frac{\partial K_i}{\partial q_1}\frac{\partial K_j}{\partial p_1}+\frac{\partial K_i}{\partial p_2}\frac{\partial K_j}{\partial q_2}-\frac{\partial K_i}{\partial q_2}\frac{\partial K_j}{\partial p_2}.
\end{equation}
Here, $\{,\}$ denotes the poisson bracket such that $\{p_i,q_j\}={\delta}_{ij}$ (${\delta}_{ij}$:kronecker's delta).

We show that this system admits the affine Weyl group symmetry of type $A_1^{(1)}$ as the group of its B{\"a}cklund transformations. We also show that its polynomial Hamiltonian system is equivalent to the following system (see Section 3):
\begin{equation}
  \left\{
  \begin{aligned}
   \frac{\partial^4 u}{\partial t_1^4} =&\frac{1}{u^2(8u^3+u-4\frac{\partial^2 u}{\partial t_1^2})}(6\left(\frac{\partial^2 u}{\partial t_1^2} \right)^2 \left(\frac{\partial u}{\partial t_1} \right)^2-2\frac{\partial^2 u}{\partial t_1^2}u \left(\frac{\partial u}{\partial t_1} \right)^2-12 \left(\frac{\partial^2 u}{\partial t_1^2} \right)^3u\\
&-4\frac{\partial^2 u}{\partial t_1^2}\frac{\partial^3 u}{\partial t_1^3}u\frac{\partial u}{\partial t_1}+7\left(\frac{\partial^2 u}{\partial t_1^2} \right)^2 u^2+2\frac{\partial^3 u}{\partial t_1^3}u^2 \frac{\partial u}{\partial t_1}-2\left(\frac{\partial^3 u}{\partial t_1^3} \right)^2 u^2-\frac{3}{2}\frac{\partial^2 u}{\partial t_1^2}u^3\\
&-48\frac{\partial^2 u}{\partial t_1^2}u^3 \left(\frac{\partial u}{\partial t_1} \right)^2+\frac{u^4}{8}+4u^4 \left(\frac{\partial u}{\partial t_1} \right)^2+24 \left(\frac{\partial^2 u}{\partial t_1^2} \right)^2 u^4+32\frac{\partial^3 u}{\partial t_1^3}u^4 \frac{\partial u}{\partial t_1}\\
&-12\frac{\partial^2 u}{\partial t_1^2}u^5+2u^6+8u^8+8\alpha_0^2 u^4),\\
   \frac{\partial u}{\partial t_2} =&-\frac{3}{2}\frac{\partial u}{\partial t_1}-\frac{\partial^3 u}{\partial t_1^3}+\frac{3\frac{\partial^2 u}{\partial t_1^2}\frac{\partial u}{\partial t_1}}{u},\\
   \frac{\partial u}{\partial t_3} =&\frac{1}{64u^3(8u^3+u-4\frac{\partial^2 u}{\partial t_1^2})}(32 \left(\frac{\partial^2 u}{\partial t_1^2} \right)^2 \left(\frac{\partial u}{\partial t_1} \right)^3-64 \left(\frac{\partial^2 u}{\partial t_1^2} \right)^3u\frac{\partial u}{\partial t_1}-64\frac{\partial^2 u}{\partial t_1^2}\frac{\partial^3 u}{\partial t_1^3}u \left(\frac{\partial u}{\partial t_1} \right)^2\\
&+48 \left(\frac{\partial^2 u}{\partial t_1^2} \right)^2u^2\frac{\partial u}{\partial t_1}+32 \left(\frac{\partial^3 u}{\partial t_1^3} \right)^2u^2\frac{\partial u}{\partial t_1}-12\frac{\partial^2 u}{\partial t_1^2}u^3\frac{\partial u}{\partial t_1}+512\frac{\partial^2 u}{\partial t_1^2}u^3 \left(\frac{\partial u}{\partial t_1} \right)^3+u^4\frac{\partial u}{\partial t_1}\\
&-64u^4 \left(\frac{\partial u}{\partial t_1} \right)^3+128 \left(\frac{\partial^2 u}{\partial t_1^2} \right)^2u^4\frac{\partial u}{\partial t_1}-256\frac{\partial^3 u}{\partial t_1^3}u^4 \left(\frac{\partial u}{\partial t_1} \right)^2-128\frac{\partial^2 u}{\partial t_1^2}u^5 \frac{\partial u}{\partial t_1}+24u^6 \frac{\partial u}{\partial t_1}\\
&+128u^8 \frac{\partial u}{\partial t_1}+32\alpha_0 \frac{\partial^2 u}{\partial t_1^2}u^4-8\alpha_0 u^5+128\alpha_0 \frac{\partial^2 u}{\partial t_1^2}u^6-96\alpha_0u^7-256\alpha_0 u^9\\
&-128\alpha_0^2 u^4\frac{\partial u}{\partial t_1}+32\alpha_1 \frac{\partial^2 u}{\partial t_1^2}u^4-8\alpha_1 u^5+128\alpha_1 \frac{\partial^2 u}{\partial t_1^2}u^6-96\alpha_1 u^7-256\alpha_1 u^9).
   \end{aligned}
  \right.
\end{equation}

\section{Symmetry and holomorphy conditions}
In this section, we study the symmetry and holomorphy conditions of the system \eqref{1}. These properties are new.

\begin{theorem}\label{th:1}
The system \eqref{1} admits the affine Weyl group symmetry of type $A_4^{(2)}$ as the group of its B{\"a}cklund transformations, whose generators $s_0,s_1,s_2$ defined as follows$:$ with {\it the notation} $(*):=(x,y,z,w,t;\alpha_0,\alpha_1,\alpha_2)$\rm{: \rm}
\begin{align}
\begin{split}
s_0:(*) \rightarrow &\left(x,y,z+\frac{\alpha_0}{w},w,t;-\alpha_0,\alpha_1+\alpha_0,\alpha_2 \right),\\
s_1:(*) \rightarrow &\left(x,y-\frac{\alpha_1}{x+z^2},z,w-\frac{2\alpha_1 z}{x+z^2},t;\alpha_0+2\alpha_1,-\alpha_1,\alpha_2+\alpha_1 \right),\\
s_2:(*) \rightarrow &\left(x+\frac{2 \alpha_2 y}{f_2}-\frac{\alpha_2^2}{f_2^2},y-\frac{\alpha_2}{f_2},z+\frac{\alpha_2}{f_2},w,t;\alpha_0,\alpha_1+2\alpha_2,-\alpha_2  \right),
\end{split}
\end{align}
where $f_2:=x+y^2+w+t$.
\end{theorem}
We note that the B{\"a}cklund transformations of this system satisfy
\begin{equation}\label{3241}
s_i(g)=g+\frac{\alpha_i}{f_i}\{f_i,g\}+\frac{1}{2!}\left(\frac{\alpha_i}{f_i} \right)^2\{f_i,\{f_i,g\} \}+\cdots \quad (g \in {\Bbb C}(t)[x,y,z,w]),
\end{equation}
where poisson bracket $\{,\}$ satisfies the relations:
$$
\{y,x\}=\{w,z\}=1, \quad the \ others \ are \ 0.
$$
Since these B{\"a}cklund transformations have Lie theoretic origin, similarity reduction of a Drinfeld-Sokolov hierarchy admits such a B{\"a}cklund symmetry.

\begin{proposition}
This system has the following invariant divisors\rm{:\rm}
\begin{center}
\begin{tabular}{|c|c|c|} \hline
parameter's relation & $f_i$ \\ \hline
$\alpha_0=0$ & $f_0:=w$  \\ \hline
$\alpha_1=0$ & $f_1:=x+z^2$  \\ \hline
$\alpha_2=0$ & $f_2:=x+y^2+w+t$  \\ \hline
\end{tabular}
\end{center}
\end{proposition}
We note that when $\alpha_0=0$, we see that the system \eqref{1} admits a particular solution $w=0$, and when $\alpha_2=0$, after we make the birational and symplectic transformations:
\begin{equation}
x_2=x+y^2+w+t, \ y_2=y, \ z_2=z+y, \ w_2=w
\end{equation}
we see that the system \eqref{1} admits a particular solution $x_2=0$.

\begin{proposition}
Let us define the following translation operators{\rm : \rm}
\begin{align}
\begin{split}
&T_1:=s_1 s_2 s_1 s_0, \quad T_2:=s_1 T_1 s_1.
\end{split}
\end{align}
These translation operators act on parameters $\alpha_i$ as follows$:$
\begin{align}
\begin{split}
T_1(\alpha_0,\alpha_1,\alpha_2)=&(\alpha_0,\alpha_1,\alpha_2)+(-2,1,0),\\
T_2(\alpha_0,\alpha_1,\alpha_2)=&(\alpha_0,\alpha_1,\alpha_2)+(0,-1,1).
\end{split}
\end{align}
\end{proposition}

\begin{theorem}\label{pro:2}
Let us consider a polynomial Hamiltonian system with Hamiltonian $K \in {\Bbb C}(t)[x,y,z,w]$. We assume that

$(A1)$ $deg(K)=6$ with respect to $x,y,z,w$.

$(A2)$ This system becomes again a polynomial Hamiltonian system in each coordinate system $r_i \ (i=0,1,2)${\rm : \rm}
\begin{align}
\begin{split}
r_0:&x_0=x, \ y_0=y, \ z_0=\frac{1}{z}, \quad w_0=-(wz+\alpha_0)z,\\
r_1:&x_1=-((x+z^2)y-\alpha_1)y, \ y_1=\frac{1}{y}, \ z_1=z, \ w_1=w-2yz,\\
r_2:&x_2=-((x+y^2+w+t)y-\alpha_2)y, \ y_2=\frac{1}{y}, \ z_2=z+y, \ w_2=w.
\end{split}
\end{align}
Then such a system coincides with the system \eqref{1} with the polynomial Hamiltonian \eqref{2}.
\end{theorem}
By this theorem, we can also recover the parameter's relation \eqref{3}.

We note that the condition $(A2)$ should be read that
\begin{align*}
&r_j(K) \quad (j=0,1), \quad r_2(K+y)
\end{align*}
are polynomials with respect to $x_i,y_i,z_i,w_i$.

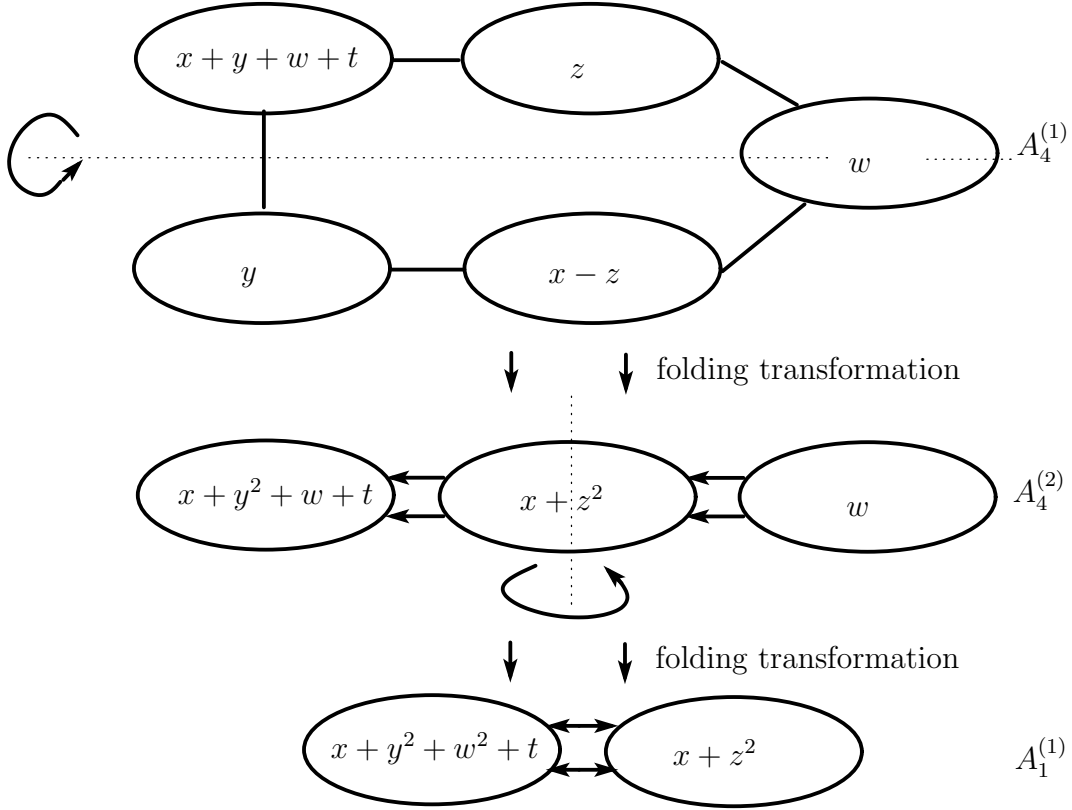
\begin{figure}[h]
\unitlength 0.1in
\begin{picture}(52.74,42.03)(3.46,-45.12)
%
\special{pn 20}%
\special{ar 1690 2882 670 286  0.0000000 6.2831853}%
%
\special{pn 20}%
\special{ar 3270 2891 670 287  0.0000000 6.2831853}%
%
\special{pn 20}%
\special{pa 2620 2791}%
\special{pa 2350 2791}%
\special{fp}%
\special{sh 1}%
\special{pa 2350 2791}%
\special{pa 2417 2811}%
\special{pa 2403 2791}%
\special{pa 2417 2771}%
\special{pa 2350 2791}%
\special{fp}%
%
\special{pn 20}%
\special{pa 2620 2993}%
\special{pa 2350 2993}%
\special{fp}%
\special{sh 1}%
\special{pa 2350 2993}%
\special{pa 2417 3013}%
\special{pa 2403 2993}%
\special{pa 2417 2973}%
\special{pa 2350 2993}%
\special{fp}%
%
\special{pn 20}%
\special{ar 4840 2892 670 287  0.0000000 6.2831853}%
%
\special{pn 20}%
\special{pa 4190 2792}%
\special{pa 3920 2792}%
\special{fp}%
\special{sh 1}%
\special{pa 3920 2792}%
\special{pa 3987 2812}%
\special{pa 3973 2792}%
\special{pa 3987 2772}%
\special{pa 3920 2792}%
\special{fp}%
%
\special{pn 20}%
\special{pa 4190 2994}%
\special{pa 3920 2994}%
\special{fp}%
\special{sh 1}%
\special{pa 3920 2994}%
\special{pa 3987 3014}%
\special{pa 3973 2994}%
\special{pa 3987 2974}%
\special{pa 3920 2994}%
\special{fp}%
\put(12.3000,-29.5900){\makebox(0,0)[lb]{$x+y^2+w+t$}}%
\put(30.2000,-29.7500){\makebox(0,0)[lb]{$x+z^2$}}%
\put(47.3000,-29.9300){\makebox(0,0)[lb]{$w$}}%
%
\special{pn 8}%
\special{pa 3290 2364}%
\special{pa 3290 3453}%
\special{dt 0.045}%
\special{pa 3290 3453}%
\special{pa 3290 3452}%
\special{dt 0.045}%
%
\special{pn 20}%
\special{pa 3102 3252}%
\special{pa 3065 3269}%
\special{pa 3030 3285}%
\special{pa 2998 3302}%
\special{pa 2970 3319}%
\special{pa 2948 3337}%
\special{pa 2934 3355}%
\special{pa 2929 3373}%
\special{pa 2933 3391}%
\special{pa 2945 3410}%
\special{pa 2965 3428}%
\special{pa 2990 3445}%
\special{pa 3022 3462}%
\special{pa 3058 3477}%
\special{pa 3098 3490}%
\special{pa 3141 3502}%
\special{pa 3187 3511}%
\special{pa 3233 3518}%
\special{pa 3281 3522}%
\special{pa 3328 3523}%
\special{pa 3374 3522}%
\special{pa 3417 3517}%
\special{pa 3458 3510}%
\special{pa 3494 3500}%
\special{pa 3526 3488}%
\special{pa 3552 3473}%
\special{pa 3570 3455}%
\special{pa 3581 3435}%
\special{pa 3584 3412}%
\special{pa 3577 3387}%
\special{pa 3562 3360}%
\special{pa 3542 3332}%
\special{pa 3532 3319}%
\special{sp}%
%
\special{pn 20}%
\special{pa 3552 3344}%
\special{pa 3482 3276}%
\special{fp}%
\special{sh 1}%
\special{pa 3482 3276}%
\special{pa 3516 3337}%
\special{pa 3520 3313}%
\special{pa 3544 3308}%
\special{pa 3482 3276}%
\special{fp}%
%
\special{pn 20}%
\special{ar 2560 4216 670 286  0.0000000 6.2831853}%
%
\special{pn 20}%
\special{ar 4140 4226 670 286  0.0000000 6.2831853}%
\put(20.3000,-42.9200){\makebox(0,0)[lb]{$x+y^2+w^2+t$}}%
\put(38.2000,-43.0900){\makebox(0,0)[lb]{$x+z^2$}}%
%
\special{pn 20}%
\special{pa 3330 4091}%
\special{pa 3190 4091}%
\special{fp}%
\special{sh 1}%
\special{pa 3190 4091}%
\special{pa 3257 4111}%
\special{pa 3243 4091}%
\special{pa 3257 4071}%
\special{pa 3190 4091}%
\special{fp}%
%
\special{pn 20}%
\special{pa 3330 4091}%
\special{pa 3500 4091}%
\special{fp}%
\special{sh 1}%
\special{pa 3500 4091}%
\special{pa 3433 4071}%
\special{pa 3447 4091}%
\special{pa 3433 4111}%
\special{pa 3500 4091}%
\special{fp}%
%
\special{pn 20}%
\special{pa 3330 4321}%
\special{pa 3190 4321}%
\special{fp}%
\special{sh 1}%
\special{pa 3190 4321}%
\special{pa 3257 4341}%
\special{pa 3243 4321}%
\special{pa 3257 4301}%
\special{pa 3190 4321}%
\special{fp}%
%
\special{pn 20}%
\special{pa 3330 4321}%
\special{pa 3500 4321}%
\special{fp}%
\special{sh 1}%
\special{pa 3500 4321}%
\special{pa 3433 4301}%
\special{pa 3447 4321}%
\special{pa 3433 4341}%
\special{pa 3500 4321}%
\special{fp}%
%
\special{pn 20}%
\special{pa 2970 3654}%
\special{pa 2970 3830}%
\special{fp}%
\special{sh 1}%
\special{pa 2970 3830}%
\special{pa 2990 3763}%
\special{pa 2970 3777}%
\special{pa 2950 3763}%
\special{pa 2970 3830}%
\special{fp}%
%
\special{pn 20}%
\special{pa 3570 3662}%
\special{pa 3570 3838}%
\special{fp}%
\special{sh 1}%
\special{pa 3570 3838}%
\special{pa 3590 3771}%
\special{pa 3570 3785}%
\special{pa 3550 3771}%
\special{pa 3570 3838}%
\special{fp}%
\put(37.3000,-38.2200){\makebox(0,0)[lb]{folding transformation}}%
\put(56.0000,-29.5900){\makebox(0,0)[lb]{$A_4^{(2)}$}}%
\put(56.1000,-43.5800){\makebox(0,0)[lb]{$A_1^{(1)}$}}%
%
\special{pn 20}%
\special{pa 2975 2137}%
\special{pa 2975 2313}%
\special{fp}%
\special{sh 1}%
\special{pa 2975 2313}%
\special{pa 2995 2246}%
\special{pa 2975 2260}%
\special{pa 2955 2246}%
\special{pa 2975 2313}%
\special{fp}%
%
\special{pn 20}%
\special{pa 3575 2146}%
\special{pa 3575 2322}%
\special{fp}%
\special{sh 1}%
\special{pa 3575 2322}%
\special{pa 3595 2255}%
\special{pa 3575 2269}%
\special{pa 3555 2255}%
\special{pa 3575 2322}%
\special{fp}%
\put(37.3500,-23.0500){\makebox(0,0)[lb]{folding transformation}}%
%
\special{pn 20}%
\special{ar 1680 596 670 287  0.0000000 6.2831853}%
\put(12.2000,-6.7300){\makebox(0,0)[lb]{$x+y+w+t$}}%
%
\special{pn 20}%
\special{ar 4850 1091 670 286  0.0000000 6.2831853}%
\put(47.4000,-11.9200){\makebox(0,0)[lb]{$w$}}%
%
\special{pn 20}%
\special{ar 3390 596 670 287  0.0000000 6.2831853}%
\put(32.8000,-6.9700){\makebox(0,0)[lb]{$z$}}%
%
\special{pn 20}%
\special{ar 1670 1694 670 287  0.0000000 6.2831853}%
\put(15.6000,-17.9400){\makebox(0,0)[lb]{$y$}}%
%
\special{pn 20}%
\special{ar 3400 1694 670 287  0.0000000 6.2831853}%
\put(31.7000,-17.8500){\makebox(0,0)[lb]{$x-z$}}%
%
\special{pn 20}%
\special{pa 4470 838}%
\special{pa 4080 613}%
\special{fp}%
%
\special{pn 20}%
\special{pa 2360 604}%
\special{pa 2700 604}%
\special{fp}%
%
\special{pn 20}%
\special{pa 1680 872}%
\special{pa 1680 1375}%
\special{fp}%
%
\special{pn 20}%
\special{pa 2350 1701}%
\special{pa 2730 1701}%
\special{fp}%
%
\special{pn 20}%
\special{pa 4080 1710}%
\special{pa 4510 1359}%
\special{fp}%
\put(56.2000,-11.4000){\makebox(0,0)[lb]{$A_4^{(1)}$}}%
%
\special{pn 8}%
\special{pa 450 1115}%
\special{pa 4630 1115}%
\special{dt 0.045}%
\special{pa 4630 1115}%
\special{pa 4629 1115}%
\special{dt 0.045}%
%
\special{pn 8}%
\special{pa 5150 1123}%
\special{pa 5590 1123}%
\special{dt 0.045}%
\special{pa 5590 1123}%
\special{pa 5589 1123}%
\special{dt 0.045}%
%
\special{pn 20}%
\special{pa 705 998}%
\special{pa 682 970}%
\special{pa 658 944}%
\special{pa 634 921}%
\special{pa 609 903}%
\special{pa 583 892}%
\special{pa 555 888}%
\special{pa 526 893}%
\special{pa 496 904}%
\special{pa 467 921}%
\special{pa 440 944}%
\special{pa 414 970}%
\special{pa 392 1000}%
\special{pa 373 1032}%
\special{pa 359 1066}%
\special{pa 350 1100}%
\special{pa 346 1133}%
\special{pa 348 1166}%
\special{pa 356 1198}%
\special{pa 370 1228}%
\special{pa 390 1255}%
\special{pa 416 1278}%
\special{pa 446 1296}%
\special{pa 478 1308}%
\special{pa 510 1312}%
\special{pa 540 1306}%
\special{pa 567 1291}%
\special{pa 591 1269}%
\special{pa 613 1243}%
\special{pa 615 1241}%
\special{sp}%
%
\special{pn 20}%
\special{pa 630 1231}%
\special{pa 710 1148}%
\special{fp}%
\special{sh 1}%
\special{pa 710 1148}%
\special{pa 649 1182}%
\special{pa 673 1186}%
\special{pa 678 1210}%
\special{pa 710 1148}%
\special{fp}%
\end{picture}%
\label{fig:A421}
\caption{The symbol in each circle denotes the invariant divisors of each system (see Theorem \ref{th:1} (resp. \ref{th:11})).}
\end{figure}

\section{Degenerate system from the system \eqref{1}}

In this section, we find a 1-parameter family of coupled Painlev\'e systems in dimension four with affine Weyl group symmetry of type $A_1^{(1)}$ given by
\begin{align}\label{11}
\begin{split}
\frac{dx}{dt}&=\frac{\partial H}{\partial y}, \quad \frac{dy}{dt}=-\frac{\partial H}{\partial x}, \quad \frac{dz}{dt}=\frac{\partial H}{\partial w}, \quad \frac{dw}{dt}=-\frac{\partial H}{\partial z}
\end{split}
\end{align}
with the polynomial Hamiltonian
\begin{align}\label{12}
\begin{split}
H =&H_{II}(x,y,t,\alpha_0)+H_3(z,w,t)+yzw\\
=&xy^2+x^2+tx-\alpha_0 y+\frac{z^2}{4}-\frac{w^2}{4}+yzw.
\end{split}
\end{align}
Here $x,y,z$ and $w$ denote unknown complex variables, and $\alpha_0,\alpha_1$ are complex parameters satisfying the relation:
\begin{equation}\label{13}
\alpha_0+\alpha_1=1
\end{equation}
The symbol $H_{3}$ is given by
\begin{align}
&H_3(z,w,t)=\frac{z^2}{4}-\frac{w^2}{4}.
\end{align}
Of course, the Hamiltonian itself is the first integral.

This is the second example which gave higher-order Painlev\'e equations of type $A_{1}^{(1)}$. 

We note that in this case the iniariant divisors are different from the ones of the second member $P_{II}^{(2)}$ of the second Painlev\'e hierarchy  given in the paper \cite{Sasa12}.
\begin{center}
\begin{tabular}{|c||c|c|} \hline 
Invariant divisors &  $f_0$ & $f_1$ \\ \hline
System \eqref{11} &  $x+z^2$ & $x+y^2+w^2+t$ \\ \hline
$P_{II}^{(2)}$ &  $y$ & $y+t-2w^2+4x(z+xw)$ \\ \hline
\end{tabular}
\end{center}

We remark that for this system we tried to seek its first integrals of polynomial type with respect to $x,y,z,w$. However, we can not find. Of course, the Hamiltonian $H$ is not its first integral.

This system can be obtained by connecting the invariant divisors $w$ and $x+y^2+w+t$ for the canonical variables $(x,y,z,w)$ in the system \eqref{1}.

\begin{theorem}\label{th:11}
The system \eqref{11} admits the affine Weyl group symmetry of type $A_1^{(1)}$ as the group of its B{\"a}cklund transformations, whose generators $s_0,s_1$ defined as follows$:$ with {\it the notation} $(*):=(x,y,z,w,t;\alpha_0,\alpha_1)$\rm{: \rm}
\begin{align}
\begin{split}
s_0:(*) \rightarrow &\left(x,y-\frac{\alpha_0}{x+z^2},z,w-\frac{2\alpha_0 z}{x+z^2},t;-\alpha_0,\alpha_1+2\alpha_0 \right),\\
s_1:(*) \rightarrow &\left(x+\frac{2 \alpha_1 y}{f_1}-\frac{\alpha_1^2}{f_1^2},y-\frac{\alpha_1}{f_1},z+\frac{2\alpha_1 w}{f_1},w,t;\alpha_0+2\alpha_1,-\alpha_1 \right),
\end{split}
\end{align}
where $f_1:=x+y^2+w^2+t$.
\end{theorem}

\begin{proposition}
This system has the following invariant divisors\rm{:\rm}
\begin{center}
\begin{tabular}{|c|c|c|} \hline
parameter's relation & $f_i$ \\ \hline
$\alpha_0=0$ & $f_0:=x+z^2$  \\ \hline
$\alpha_1=0$ & $f_1:=x+y^2+w^2+t$  \\ \hline
\end{tabular}
\end{center}
\end{proposition}

\begin{proposition}
Let us define the following translation operator{\rm : \rm}
\begin{align}
\begin{split}
&T:=s_1 s_0.
\end{split}
\end{align}
This translation operator acts on parameters $\alpha_i$ as follows$:$
\begin{align}
\begin{split}
T(\alpha_0,\alpha_1)=&(\alpha_0,\alpha_1)+(-2,2).
\end{split}
\end{align}
\end{proposition}

\begin{theorem}\label{pro:22}
Let us consider a polynomial Hamiltonian system with Hamiltonian $K \in {\Bbb C}(t)[x,y,z,w]$. We assume that

$(B1)$ $deg(K)=6$ with respect to $x,y,z,w$.

$(B2)$ This system becomes again a polynomial Hamiltonian system in each coordinate system $r_i \ (i=0,1)${\rm : \rm}
\begin{align}
\begin{split}
r_0:&x_0=-((x+z^2)y-\alpha_0)y, \ y_0=\frac{1}{y}, \ z_0=z, \quad w_0=w-2yz,\\
r_1:&x_1=-((x+y^2+w^2+t)y-\alpha_1)y, \ y_1=\frac{1}{y}, \ z_1=z+2yw, \ w_1=w.
\end{split}
\end{align}
Then such a system coincides with the system \eqref{11} with the polynomial Hamiltonian \eqref{12}.
\end{theorem}
By this theorem, we can also recover the parameter's relation \eqref{13}.

We note that the condition $(B2)$ should be read that
\begin{align*}
&r_0(K), \quad r_1(K+y)
\end{align*}
are polynomials with respect to $x_i,y_i,z_i,w_i$.

\section{1-parameter family of partial differential systems in three variables}

In this section, we find a 1-parameter family of partial differential systems in three variables given by
\begin{equation}\label{eq:10}
  \left\{
  \begin{aligned}
   dq_1 =&\frac{\partial K_1}{\partial p_1}dt_1+\frac{\partial K_2}{\partial p_1}dt_2+\frac{\partial K_3}{\partial p_1}dt_3,\\
   dp_1 =&-\frac{\partial K_1}{\partial q_1}dt_1-\frac{\partial K_2}{\partial q_1}dt_2-\frac{\partial K_3}{\partial q_1}dt_3,\\
   dq_2 =&\frac{\partial K_1}{\partial p_2}dt_1+\frac{\partial K_2}{\partial p_2}dt_2+\frac{\partial K_3}{\partial p_2}dt_3,\\
   dp_2 =&-\frac{\partial K_1}{\partial q_2}dt_1-\frac{\partial K_2}{\partial q_2}dt_2-\frac{\partial K_3}{\partial q_2}dt_3
   \end{aligned}
  \right. 
\end{equation}
with the polynomial Hamiltonians
\begin{align}
\begin{split}
K_1=&q_1p_1^2+q_1^2-\alpha_0 p_1+\frac{1}{4}q_2^2-\frac{1}{4}p_2^2+p_1q_2p_2,\\
K_2=&q_2^2p_2^2-\frac{1}{4}q_2^2+\frac{1}{4}p_2^2-2\alpha_0 q_2p_2+q_1q_2^2+q_1p_2^2-p_1q_2p_2+p_1^2q_2^2,\\
K_3=&\frac{1}{2}q_1^2p_1^4+q_1^3p_1^2+\frac{1}{2}q_1^4-\alpha_0 q_1p_1^3+\alpha_1 q_1^2p_1+\frac{1}{2}\alpha_0^2 p_1^2+\frac{1}{32}q_2^4+\frac{1}{32}p_2^4-\frac{1}{16}q_2^2p_2^2+q_1p_1^3q_2p_2\\
&+\frac{1}{2}p_1^2q_2^2p_2^2-\frac{1}{4}q_1p_1^2p_2^2-\frac{1}{4}p_1q_2p_2^3+q_1^2p_1q_2p_2+\frac{1}{4}q_1p_1^2q_2^2+\frac{1}{4}p_1q_2^3p_2-\alpha_0 p_1^2q_2p_2\\
&+\frac{1}{4}\alpha_0 p_1p_2^2+\frac{1}{4}\alpha_1 p_1q_2^2.
\end{split}
\end{align}

\begin{proposition}
The system \eqref{eq:10} satisfies the compatibility conditions$:$
\begin{equation}
\frac{\partial }{\partial t_j} \frac{\partial q_1}{\partial t_i}=\frac{\partial }{\partial t_i} \frac{\partial q_1}{\partial t_j}, \quad \frac{\partial }{\partial t_j} \frac{\partial p_1}{\partial t_i}=\frac{\partial }{\partial t_i} \frac{\partial p_1}{\partial t_j}, \quad \frac{\partial }{\partial t_j} \frac{\partial q_2}{\partial t_i}=\frac{\partial }{\partial t_i} \frac{\partial q_2}{\partial t_j}, \quad \frac{\partial }{\partial t_j} \frac{\partial p_2}{\partial t_i}=\frac{\partial }{\partial t_i} \frac{\partial p_2}{\partial t_j},
\end{equation}
where $i,j=1,2,3$.
\end{proposition}

\begin{proposition}
The system \eqref{eq:10} has $K_1,K_2$ and $K_3$ as its first integrals.
\end{proposition}

\begin{proposition}
Three Hamiltonians $K_1,K_2$ and $K_3$ satisfy
\begin{equation}
\{K_1,K_2\}=\{K_1,K_3\}=\{K_2,K_3\}=0,
\end{equation}
where
\begin{equation}
\{K_i,K_j\}=\frac{\partial K_i}{\partial p_1}\frac{\partial K_j}{\partial q_1}-\frac{\partial K_i}{\partial q_1}\frac{\partial K_j}{\partial p_1}+\frac{\partial K_i}{\partial p_2}\frac{\partial K_j}{\partial q_2}-\frac{\partial K_i}{\partial q_2}\frac{\partial K_j}{\partial p_2}.
\end{equation}
\end{proposition}
Here, $\{,\}$ denotes the poisson bracket such that $\{p_i,q_j\}={\delta}_{ij}$ (${\delta}_{ij}$:kronecker's delta).

\begin{theorem}
The system \eqref{eq:10} admits the affine Weyl group symmetry of type $A_1^{(1)}$ as the group of its B{\"a}cklund transformations, whose generators $s_0,s_1$ defined as follows$:$ with {\it the notation} $(*):=(q_1,p_1,q_2,p_2,t_1,t_2,t_3;\alpha_0,\alpha_1)$\rm{: \rm}
\begin{align}
\begin{split}
s_0:(*) \rightarrow &\left(q_1,p_1-\frac{\alpha_0}{q_1+q_2^2},q_2,p_2-\frac{2\alpha_0 q_2}{q_1+q_2^2},t_1,t_2,t_3;-\alpha_0,\alpha_1+2\alpha_0 \right),\\
s_1:(*) \rightarrow &\left(q_1+\frac{2 \alpha_1 p_1}{f_1}-\frac{\alpha_1^2}{f_1^2},p_1-\frac{\alpha_1}{f_1},q_2+\frac{2\alpha_1 p_2}{f_1},p_2,t_1,t_2,t_3;\alpha_0+2\alpha_1,-\alpha_1 \right),
\end{split}
\end{align}
where $f_1:=q_1+p_1^2+p_2^2$.
\end{theorem}
Here, the parameters $\alpha_i$ satisfy the relation $\alpha_0+\alpha_1=0$.

\begin{theorem}
Let us consider a polynomial Hamiltonian system with Hamiltonian $K \in {\Bbb C}[q_1,p_1,q_2,p_2]$. We assume that

$(C1)$ $deg(K)=6$ with respect to $q_1,p_1,q_2,p_2$.

$(C2)$ This system becomes again a polynomial Hamiltonian system in each coordinate $R_i \ (i=0,1)${\rm : \rm}
\begin{align*}
\begin{split}
R_0:(x_0,y_0,z_0,w_0)=&\left(-((q_1+q_2^2)p_1-\alpha_0)p_1,\frac{1}{p_1},q_2,p_2-2p_1q_2 \right),\\
R_1:(x_1,y_1,z_1,w_1)=&\left(-((q_1+p_1^2+p_2^2)p_1-\alpha_1)p_1,\frac{1}{p_1},q_2+2p_1p_2,p_2 \right),
\end{split}
\end{align*}
where the parameters $\alpha_i$ satisfy the relation $\alpha_0+\alpha_1=0$. Then such a system coincides with the Hamiltonian system \eqref{eq:10} with the polynomial Hamiltonians $K_1,K_2$ and $K_3$.
\end{theorem}
We note that the conditions $(C2)$ should be read that
\begin{align*}
\begin{split}
&R_0(K), \quad R_1(K)
\end{split}
\end{align*}
are polynomials with respect to $x_i,y_i,z_i,w_i$.

Next, let us consider the relation between the polynomial Hamiltonian system \eqref{eq:10} and and an autonomous version of the system \eqref{11}. In this paper, we can make the birational transformations between the polynomial Hamiltonian system \eqref{eq:10} and  an autonomous version of the system \eqref{11}.
\begin{theorem}
The birational transformations
\begin{equation}\label{eq:11}
  \left\{
  \begin{aligned}
   x =&\frac{1}{4}q_2-2q_1q_2,\\
   y =&-\frac{1}{8}p_2+q_1p_2+\frac{1}{4}p_1q_2-6q_1p_1q_2-2q_2^2p_2+2\alpha_0 q_2,\\
   z =&q_2,\\
   w =&-\frac{1}{2}p_2+p_1q_2
   \end{aligned}
  \right. 
\end{equation}
take the Hamiltonian system \eqref{eq:10} to the system
\begin{equation}\label{eq:12}
  \left\{
  \begin{aligned}
   dx =&y dt_1+f_1(x,y,z,w)dt_2+f_2(x,y,z,w)dt_3,\\
   dy =&\frac{1}{z^2(8z^3+z-4x)}(6x^2w^2-2xzw^2-12x^3z-4xyzw+7x^2 z^2+2yz^2 w-2y^2 z^2-\frac{3}{2}xz^3\\
&-48xz^3w^2+\frac{z^4}{8}+4z^4 w^2+24x^2 z^4+32yz^4 w-12xz^5+2z^6+8z^8+8\alpha_0^2 z^4) dt_1\\
&+\left(-\frac{3}{2}w-y+\frac{3xw}{z} \right)dt_2\\
&+\frac{1}{64z^3(8z^3+z-4x)}(32x^2w^3-64x^3zw-64xyzw^2+48x^2z^2w+32y^2z^2w-12xz^3w\\
&+512xz^3w^3+z^4w-64z^4w^3+128x^2z^4w-256yz^4w^2-128xz^5 w+24z^6 w+128z^8 w\\
&+32\alpha_0 xz^4-8\alpha_0 z^5+128\alpha_0 xz^6-96\alpha_0z^7-256\alpha_0 z^9-128\alpha_0^2 z^4w+32\alpha_1 xz^4-8\alpha_1 z^5\\
&+128\alpha_1 xz^6-96\alpha_1 z^7-256\alpha_1 z^9)dt_3,\\
   dz =&w dt_1+f_3(x,y,z,w)dt_2+f_4(x,y,z,w)dt_3,\\
   dw =&x  dt_1+f_5(x,y,z,w)dt_2+f_6(x,y,z,w)dt_3,
   \end{aligned}
  \right. 
\end{equation}
where $f_i(x,y,z,w) \in {\Bbb C}[x,y,z,w]$.
\end{theorem}
Setting $u:=z$, we see that
\begin{equation}
\frac{\partial u}{\partial t_1}=w, \quad  \frac{\partial^2 u}{\partial t_1^2}=x, \quad \frac{\partial^3 u}{\partial t_1^3}=y,
\end{equation}
and
\begin{equation}\label{eq:13}
  \left\{
  \begin{aligned}
   \frac{\partial^4 u}{\partial t_1^4} =&\frac{1}{u^2(8u^3+u-4\frac{\partial^2 u}{\partial t_1^2})}(6\left(\frac{\partial^2 u}{\partial t_1^2} \right)^2 \left(\frac{\partial u}{\partial t_1} \right)^2-2\frac{\partial^2 u}{\partial t_1^2}u \left(\frac{\partial u}{\partial t_1} \right)^2-12 \left(\frac{\partial^2 u}{\partial t_1^2} \right)^3u\\
&-4\frac{\partial^2 u}{\partial t_1^2}\frac{\partial^3 u}{\partial t_1^3}u\frac{\partial u}{\partial t_1}+7\left(\frac{\partial^2 u}{\partial t_1^2} \right)^2 u^2+2\frac{\partial^3 u}{\partial t_1^3}u^2 \frac{\partial u}{\partial t_1}-2\left(\frac{\partial^3 u}{\partial t_1^3} \right)^2 u^2-\frac{3}{2}\frac{\partial^2 u}{\partial t_1^2}u^3\\
&-48\frac{\partial^2 u}{\partial t_1^2}u^3 \left(\frac{\partial u}{\partial t_1} \right)^2+\frac{u^4}{8}+4u^4 \left(\frac{\partial u}{\partial t_1} \right)^2+24 \left(\frac{\partial^2 u}{\partial t_1^2} \right)^2 u^4+32\frac{\partial^3 u}{\partial t_1^3}u^4 \frac{\partial u}{\partial t_1}\\
&-12\frac{\partial^2 u}{\partial t_1^2}u^5+2u^6+8u^8+8\alpha_0^2 u^4),\\
   \frac{\partial u}{\partial t_2} =&-\frac{3}{2}\frac{\partial u}{\partial t_1}-\frac{\partial^3 u}{\partial t_1^3}+\frac{3\frac{\partial^2 u}{\partial t_1^2}\frac{\partial u}{\partial t_1}}{u},\\
   \frac{\partial u}{\partial t_3} =&\frac{1}{64u^3(8u^3+u-4\frac{\partial^2 u}{\partial t_1^2})}(32 \left(\frac{\partial^2 u}{\partial t_1^2} \right)^2 \left(\frac{\partial u}{\partial t_1} \right)^3-64 \left(\frac{\partial^2 u}{\partial t_1^2} \right)^3u\frac{\partial u}{\partial t_1}-64\frac{\partial^2 u}{\partial t_1^2}\frac{\partial^3 u}{\partial t_1^3}u \left(\frac{\partial u}{\partial t_1} \right)^2\\
&+48 \left(\frac{\partial^2 u}{\partial t_1^2} \right)^2u^2\frac{\partial u}{\partial t_1}+32 \left(\frac{\partial^3 u}{\partial t_1^3} \right)^2u^2\frac{\partial u}{\partial t_1}-12\frac{\partial^2 u}{\partial t_1^2}u^3\frac{\partial u}{\partial t_1}+512\frac{\partial^2 u}{\partial t_1^2}u^3 \left(\frac{\partial u}{\partial t_1} \right)^3+u^4\frac{\partial u}{\partial t_1}\\
&-64u^4 \left(\frac{\partial u}{\partial t_1} \right)^3+128 \left(\frac{\partial^2 u}{\partial t_1^2} \right)^2u^4\frac{\partial u}{\partial t_1}-256\frac{\partial^3 u}{\partial t_1^3}u^4 \left(\frac{\partial u}{\partial t_1} \right)^2-128\frac{\partial^2 u}{\partial t_1^2}u^5 \frac{\partial u}{\partial t_1}+24u^6 \frac{\partial u}{\partial t_1}\\
&+128u^8 \frac{\partial u}{\partial t_1}+32\alpha_0 \frac{\partial^2 u}{\partial t_1^2}u^4-8\alpha_0 u^5+128\alpha_0 \frac{\partial^2 u}{\partial t_1^2}u^6-96\alpha_0u^7-256\alpha_0 u^9\\
&-128\alpha_0^2 u^4\frac{\partial u}{\partial t_1}+32\alpha_1 \frac{\partial^2 u}{\partial t_1^2}u^4-8\alpha_1 u^5+128\alpha_1 \frac{\partial^2 u}{\partial t_1^2}u^6-96\alpha_1 u^7-256\alpha_1 u^9).
   \end{aligned}
  \right.
\end{equation}

\end{document}